\documentclass[A4paper, 12pt]{article}
\usepackage{latexsym}
\usepackage{indentfirst}
\usepackage{listings}
\usepackage{caption}
\usepackage{mathtext}
\usepackage[dvips]{graphicx}
\usepackage{epsfig}
\usepackage{amssymb,amsthm,amsmath,amsfonts, amscd}

\usepackage[top=3cm, bottom=3cm, left=3cm, right=2cm]{geometry}
\makeatletter
\renewcommand\@biblabel[1]{#1.}
\makeatother
\long\def\@makecaption#1#2{%
  \vskip\abovecaptionskip
  \sbox\@tempboxa{#1.~#2}%
  \ifdim \wd\@tempboxa >\hsize
    #1.~#2\par
  \else
    \global \@minipagefalse
    \hb@xt@\hsize{\hfil\box\@tempboxa\hfil}%
  \fi
  \vskip\belowcaptionskip}
\reversemarginpar

\begin{document}

\noindent
{\Large \bf  New $U$-empirical tests of symmetry based on extremal order statistics, and their efficiencies}

\bigskip

\noindent Ya.~Yu.~NIKITIN$^{a,}\footnote[1]{Supported by RFBR grant No. 13-01-00172,
and by SPbGU grant No. 6.38.672.2013.}$ and M.~AHSANULLAH$^{b}$ \\

\medskip

\noindent {\it
$^{a}$Department of Mathematics and Mechanics, Saint-Petersburg State University, Universitetsky pr. 28,
Stary Peterhof 198504, Russia, and National Research University - Higher School of Economics, Souza
Pechatnikov, 16, St.Peters\-burg 190008, Russia;  \\}

\vskip1pt

\noindent{\it  $^{b}$Department of Management Sciences, Rider University, Lawrenceville, NJ 08648, USA.}

\bigskip

\noindent E-mail:  yanikit47@gmail.com   \, , \, ahsan@rider.edu

\medskip

\noindent Running title: New tests of symmetry.

\bigskip

\noindent
{\it \small Abstract. We use a characterization of symmetry in terms of extremal order statistics which enables to build several new nonparametric tests of symmetry. We discuss their limiting distributions and calculate their local exact Bahadur efficiency under location alternative which is mostly high.}

\medskip


\medskip


\medskip

\section{Introduction}

\noindent The idea of building statistical tests based on characterizations belongs to Yu.V.Linnik \cite{Lin}. Suppose we have a sample  $X_1,\ldots,X_n$ of i.i.d. observations with the df $F,$ and we are testing the hypothesis  ${\cal H}: F \in {\cal F},$ where $\cal F$ is some family of distributions,  against the alternative ${\cal A} : F\notin {\cal F}.$ Common examples of $\cal F$ are the families of exponential or normal distributions with unknown parameters or the class of symmetric distributions with known or unknown center of symmetry.

Consider the characterization of $\cal F$ by equal distribution of two statistics  $g_1(X_1, \ldots , X_r)$ and $g_2(X_1, \ldots , X_s).$
We introduce  two  $U$-empirical df's
\begin{gather*}
G_{1n}(t)={n \choose r}^{-1}\sum_{1 \leq i_1<\ldots <i_r \leq n}\textbf{1}\{g_1(X_{i_1}, \ldots , X_{i_r})< t\}, \quad t\in R^1,\quad r \geq 1,\\
G_{2n}(t)={n \choose s}^{-1}\sum_{1 \leq i_1<\ldots <i_s \leq n}\textbf{1}\{g_2(X_{i_1}, \ldots , X_{i_s})<t\},\quad t\in R^1, \quad s \geq 1.
\end{gather*}

According to the Glivenko-Cantelli theorem for $U$-empirical df's, see \cite{helm}, $G_{jn}(t)$ converge uniformly and a.s. to the df's $G_j (t)=P (g_j <t ),  j=1,2, $ as $n \to \infty.$ As under ${\cal H}$ one has  $G_1(t)\equiv G_2 (t),$ it follows that a.s. under ${\cal H}$
$$
D_n := \sup_{t\in R^1}\mid G_{1n}(t)- G_{2n}(t)\mid \longrightarrow 0, \,\, n\to \infty.
$$
Hence the Kolmogorov-type statistic $D_n$ can be used for testing ${\cal H}$  against ${\cal A}.$ We can also use some $U$-empirical integral statistics, e.g.
\begin{equation}
\label{integral}
I_n=\int_{R} \left(G_{1n}(t)-G_{2n}(t)\right)dF_n(t),
\end{equation}
where $F_n $ is the usual empirical df, {\it in case they are consistent.} The use of $\omega^2$-type statistics of the type $$
\Omega_n^2=\int_{R} \left(G_{1n}(t)-G_{2n}(t)\right)^2 dF_n(t) $$
is likely to be unjustified because of their complexity and considerable difficulty of applying limit theory.

The examples of such  {\it goodness-of-fit tests} together with their asymptotic analysis and related calculation of efficiencies can be found in \cite{bih}, \cite{henMei}, \cite{morSz},  \cite{mulnik}, \cite{nik96},  \cite{nikolm}, \cite{nikVol}, and some other related papers.

Testing of {\it symmetry} based on characterizations has been much less explored. Consider the classical hypothesis
\begin{equation}
\label{symm}
H_0: 1 - F(x) - F(-x) = 0, \, \, \forall x \in R^1,
\end{equation}
against the alternative $H_1$ under which the equality (\ref{symm}) is violated at least in one point. The first step in construction of such tests was made by Baringhaus and Henze in \cite{bih}.

Suppose that $X$ and $Y$ are i.i.d.  rv's with continuous df $F$. Baringhaus and Henze proved that the equal distribution of rv's $|X|$ and $|\max(X,Y)|$ is valid iff $F$ is symmetric with respect to zero, that is (\ref{symm}) holds. They also proposed suitable Kolmogorov-type  and omega-square type tests of symmetry. Some efficiency calculations were then performed in \cite{nikMMS}, see also \cite{nikolm}. Integral test of symmetry similar to (\ref{integral}) was proposed next by Litvinova \cite{lit}. In the present paper we reconsider { \it inter alia} the Litvinova test.

In our paper we are interested  in {\it new tests of symmetry} with respect to zero based on the following characterization by Ahsanullah \cite{ahsan}:

\, \textit{Suppose that $X_{1},...,X_{k}, k \geq 2,$ are i.i.d. rv's with absolutely continuous df $F(x).$  Denote $X_{1,k}= \min (X_{1},...,X_{k})$ and $X_{k,k}=\max (X_{1},...,X_{k}).$ Then  $|X_{1,k}|$ and $|X_{k,k}|$ are identically distributed iff $F$ is symmetric about zero, i.e. (\ref{symm}) holds.}

Subsequently we refer to this result as {\it Ahsanullah's characterization of order $k$.}

In the sequel we construct new tests of symmetry using this characterization and explore their asymptotic properties with emphasis on their local Bahadur efficiency.  We shall see that corresponding tests of symmetry for $ k =2$ and $k = 3$ are asymptotically equivalent to the test of Litvinova and to the Kolmogorov-type test of Baringhaus and Henze. In case of location alternative they are competitive and manifest rather high Bahadur and Pitman efficiency in comparison to many other tests of symmetry. At the same time, higher values of $k, k >3,$ lead us to different tests with presumably lower values of efficiencies in case of common alternatives.

In the rest of the Introduction we present brief information on asymptotic normality of $U$-statistics and a short background on the calculation of Bahadur efficiency which is repeatedly used later on and might be helpful for the reader.

Currently $U$-statistics play an important role in Statistics and Probability. They appeared in the middle of 1940-s in problems of unbiased estimation \cite{halm}, but after the seminal paper of Hoeffding \cite{hoe} it became clear that the numerous valuable statistics are just $U$-statistics (or von Mises functionals having very similar asymptotic theory.) Most complete exposition of theory can be found in monographs \cite{kor} and \cite{Lee}.

We consider $U-$statistics of the form
$$
    U_n={n \choose m}^{-1}\sum_{1\leqslant i_1<\ldots<i_m\leqslant
    n}{\Psi(X_{i_1},\ldots,X_{i_m})},\qquad n\geqslant m ,
$$
where $X_1,X_2,\dots$ is a sequence of i.i.d. rv's with common distribution $P$, while the kernel $\Psi:R^m \to R^1$ is a measurable symmetric function of $m$ variables. The number $m$ is called the {\it degree} of the kernel. We assume that the kernel $\Psi$ is integrable on $R^m$ and denote
$$
\theta(P) = \int...\int_{R^m} \Psi(x_1,\ldots, x_m) dP(x_1)... \ dP(x_m).
$$
In the sequel we need the notations
$$
    \psi(x):= \mathbb{E}_P\{\Psi(X_1,\ldots,X_m)|X_1 = x\}, \quad
    \Delta^2 := \mathbb{E}_P \psi^2(X_1) - (\theta(P))^2.
    $$

The function  $\psi$ is called the one-dimensional {\it projection} of the kernel $\Psi$ and plays an important role in asymptotic theory. If $\Delta^2 > 0$ that specifies the so-called non-degenerate case, the limiting distribution of $U-$statistics is normal as discovered by Hoeffding \cite{hoe}. He proved that
if $ \mathbb{E}_P \Psi^2(X_1,\ldots,X_m) < \infty$ and $\Delta^2 > 0,$ then as $n \to \infty$ one has convergence in distribution
 \begin{equation}
 \label{normal}
    \sqrt{\frac{n}{m^2\Delta^2}}\left(U_n - \theta(P)\right)\mathrel{\stackrel{\makebox[0pt]{\mbox{\normalfont\tiny d}}}{\longrightarrow}} N(0,1).
\end{equation}

Bahadur efficiency is one of several possible approaches to measure the
asymptotic relative efficiency (ARE) of statistical tests.
The Bahadur approach proposed in \cite{Bah67}, \cite{bah}
 prescribes one to fix the power of concurrent tests and to
compare the exponential rates of decrease of their sizes for  the
increasing number of observations and fixed alternative. This
exponential rate for a sequence of statistics $\{T_n\}$ is usually
proportional to some non-random function $c_T(\theta)$ depending
on the alternative parameter $\theta$ which is called the {\it
exact slope} of the sequence $\{T_n\}$. The Bahadur ARE $\,
e_{V,T}^{\, B} (\theta)$ of two sequences of statistics
$\,\{V_n\}$ and $\,\{T_n\}$ is defined by means of the formula
$$e_{V,T}^{\,B}(\theta) = c_V(\theta)\,\big/\,c_T(\theta)\,.$$

The Bahadur exact slope of the sequence of test statistics $\{T_n\}$ can be evaluated as
$c_T (\theta) = 2f(b_T(\theta)),$ where $b_T(\theta)$ is the limit in probability of $T_n$
under the alternative, while the continuous function $f(t)$ describes the logarithmic large deviation
asymptotics of this sequence under the null-hypothesis, see details in \cite{bah} or \cite{nik}.

\par  It is important to note that there exists an upper bound  for  exact  slopes \cite{Bah67}, \cite{bah}
\begin{equation}
\label{KL}
c_T(\theta) \leq 2K(\theta),
\end{equation}
where the Kullback--Leibler information number $K(\theta)$ measures the "statistical
distance" between the alternative and the null-hypothesis. It is
sometimes compared in the literature with the Cram\'er--Rao
inequality in the estimation theory. Therefore the absolute
(nonrelative) Bahadur efficiency of the sequence $\{T_n\}$  can be
defined as $e_T^B(\theta) := c_T(\theta)/2K(\theta).$

Often the exact Bahadur ARE is uncomputable for any alternative depending on
$\theta$ but it is possible to calculate the local Bahadur ARE as
$\theta$ approaches the null-hypothesis. Then one speaks about the
{\it local} Bahadur efficiency and {\it local} Bahadur slopes \cite{nik}.

The indisputable merit of Bahadur efficiency in the ability to be calculated for statistics with non-normal asymptotic
distribution. This is the primary reason to use it in the present paper as the Kolmogorov-type statistics
have non-normal limiting distribution.

\section{Integral test of symmetry for $k=2$ and its asymptotic theory}

\noindent In this section we study the simplest integral test. Consider two $V$-empirical df's
$$
\begin{array}{ll}
\medskip
G_n(t) = n^{-2} \sum_{1\le i,j\le n} {\bf{1}}\{ |\min(X_i, X_j)| < t\}, t\in R^1, \\
H_n(t) = n^{-2} \sum_{1\le i,j\le n} {\bf{1}}\{ |\max(X_i, X_j)| < t\}, t\in R^1,
\end{array}
$$
and let $Q_n$ be the empirical df corresponding to the sample $|X_i|, i =1,...,n.$

We introduce the integral statistic
$$
J_n = \int_{R^1} [G_n(t)- H_n(t)]dQ_n(t).
$$
Now we shall show that this statistic is distribution-free under the hypothesis of symmetry. Denote by $F^{-1}$ the inverse df \ of the sample assuming for simplicity that it is strictly monotone. Then
$$
J_n = \int_0^1 [H_n( F^{-1} (u)) - G_n(F^{-1} (u))]dQ_n(F^{-1} (u)).
$$
By symmetry of general df $F$,
$$- F^{-1} (u) = F^{-1} (1-u), \qquad \forall u \in [0,1],$$
hence for any $u$
\begin{multline*}
G_n( F^{-1} (u)) = n^{-2} \sum_{1\le i,j\le n} {\bf{1}}\{ - F^{-1} (u) <\min(X_i, X_j) < F^{-1} (u)\}  \\ =
n^{-2} \sum_{1\le i,j\le n} {\bf{1}}\{  F^{-1} (1-u) <\min(X_i, X_j) < F^{-1} (u)\} \\= n^{-2} \sum_{1\le i,j\le n} {\bf{1}}\{  (1-u) <\min(F(X_i),F(X_j)) < u\} \\= n^{-2} \sum_{1\le i,j\le n} {\bf{1}}\{  1-u <\min(U_i, U_j) < u\},
\end{multline*}
where $U_1,...,U_n$ are independent standard uniform rv's, and we see that  $G_n( F^{-1} (u))$ does not depend on $F$. Similar arguments are true for
$H_n(F^{-1} (u))$ and $Q_n(F^{-1} (u)).$ Hence $J_n$ is distribution-free. Thus we may assume in the sequel that $F$ is the (symmetric) uniform distribution on $[-1,1].$ Now we see that
\begin{multline*}
J_n = n^{-3} \sum_{1\le i,j,k\le n} \left({\bf{1}}\{ |\min(X_i, X_j)| < |X_k|\} -  {\bf{1}}\{ |\max(X_i, X_j)| < |X_k|\}\right)\\= n^{-3} \sum_{1\le i,j,k\le n} \Psi_3(X_i, X_j, X_k),
 \end{multline*}
where the kernel $\Psi_3$ of degree 3 of the last $V$-statistic  is given after symmetrization by
 \begin{multline*}
 3\Psi_3 (X, Y, Z) = {\bf{1}}\{ |\min(X,Y)| < |Z|\} + {\bf{1}}\{ |\min(X,Z)| < |Y|\}+ {\bf{1}}\{ |\min(Y,Z)| < |X|\}  \\- {\bf{1}}\{ |\max(X,Y)|< |Z|\}-{\bf{1}}\{ |\max(X,Z)|< |Y|\}- {\bf{1}}\{ |\max(Y,Z)|< |X|\}.
 \end{multline*}

 As $U-$ and $V-$statistics with the same kernel have the same asymptotic distribution \cite{kor}, we can replace the $V-$statistic $J_n$ by asymptotically equivalent $U-$statistic $I_n^{(3)}$ of degree 3
 $$
 I_n^{(3)} = { n\choose 3 }^{-1}\sum_{1\le i< j < k \le n} \Psi_3(X_i, X_j, X_k),
 $$
which is simpler to calculate.

In what follows we use a system of notations for statistics $I_n^{(k)}$ and $D_n^{(k)}$ in such a way that the index $k$ always corresponds to the {\it degree} of associated $U$-statistic or to the {\it degree} of corresponding {\it family} of $U$-statistics in case of supremum type tests. At the same time these statistics correspond to the Ahsanullah's characterization {\it of order $k-1$.}

 Let us calculate the projection of the kernel $\Psi_3.$  We should find
 $$
 \psi_3 (s) := {\mathbb E}[ \Psi_3 (X, Y, Z)| Z=s].
  $$
 Due to the underlying characterization, we have
 $$
 {\mathbb E} \left( {\bf{1}}\{ |\min(X,Y)| < |s|\} - {\bf{1}}\{ |\max(X,Y)| < |s|\} \right)  =0.
 $$
 It is clear that
 $$
 {\mathbb E} \ {\bf{1}}\{ |\min(X,s)| < |Y|\} = {\mathbb E} \ {\bf{1}}\{ |\min(Y,s)| < |X|\} =
 {\mathbb P}\{ |\min(X,s)| < |Y|\}.
 $$
The simplest way to calculate this probability is to use geometric considerations, evaluating
\begin{multline*}
\medskip
 \frac14 meas\{(x,y): -1\le x, y \le 1, |\min(x,s)| < |y|\} = \left\{ \begin{array} {lll} (s^2 -2s +3)/4, &\mbox {\rm if} & s> 0
,\\ (-s^2 +2s +3)/4, &\mbox {\rm if} &s \leq 0.
\end{array}\right.
\end{multline*}

The values of the expectations
$$
{\mathbb E}{\bf{1}}\{ |\max(X,s)|< |Y|\} = {\mathbb E}{\bf{1}}\{ |\max(Y,s)|< |X|\} = {\mathbb P}\{ |\max(X,s)|< |Y|\}
$$
are slightly different and are given by
\begin{multline*}
\medskip
 \frac14 meas\{(x,y): -1\le x, y \le 1, |\max(x,s)| < |y|\} = \left\{ \begin{array} {lll} (-s^2 -2s +3)/4, &\mbox {\rm if} &s> 0
,\\ (s^2 +2s +3)/4, &\mbox {\rm if} &s \leq 0.
\end{array}\right.
\end{multline*}

Hence
$$
{\mathbb E}{\bf{1}}\{ |\min(X,s)|< |Y|\} - {\mathbb E}{\bf{1}}\{ |\max(X,s)|< |Y|\} = \left\{ \begin{array} {lll} &s^2 /2, \, \mbox {\rm if} & s> 0
,\\ -&s^2 /2,  \, \mbox {\rm if} &s \leq 0.
\end{array}\right.
$$
Taking in account the same value for ${\mathbb E}  {\bf{1}}\{ |\min(Y,s)| < |X|\} -{\mathbb E}  {\bf{1}}\{ |\max(Y,s)| < |X|\}$, we conclude that the required projection is given by
$$
 \psi_3 (s)=\left\{ \begin{array} {lll} &s^2 /3, \,\, \mbox {\rm if} & s> 0
,\\ -&s^2 /3,  \,\, \mbox {\rm if} &s \leq 0.
\end{array}\right.
$$

Consequently the projection's variance equals
$$
\sigma_3^2 := {\mathbb E} \psi_3^2 (X_1) = \frac{1}{18} \int_{-1}^{1} x^4 \ dx = \frac{1}{45} > 0,
$$
so that our kernel $\Psi_3$ is non-degenerate. According to Hoeffding's theorem, see (\ref{normal}), we have weak convergence
$$
\sqrt{n}I_n^{(3)} \mathrel{\stackrel{\makebox[0pt]{\mbox{\normalfont\tiny d}}}{\longrightarrow}} N (0, \ \frac{1}{5} ).
$$

Now we can  describe the rough large deviation asymptotics under $H_0.$ The following result is proved in \cite{nikiponi}: \textit{
 For $a>0$ it holds true under $H_0$ that
$$
\lim_{n\to \infty} n^{-1} \ln P( I_n^{(3)} >a) = - f_3 (a),
$$
where the function $f_3$ is analytic for sufficiently small $a>0,$ and that $$
f_3 (a) \sim \frac{a^2}{18 \sigma^2} = \frac52 \, a^2, \quad \mbox{as} \,
\, a \to 0.
$$}

Now we apply Bahadur's theory \cite{bah}, \cite{nik}  to evaluate the local Bahadur efficiency of this test. By the Law of Large Numbers for $U$- and $V$-statistics, see \cite{kor}, we have a.s. convergence under the parametric alternative $P_{\theta}$:
$$
I_n^{(3)} \mathrel{\stackrel{\makebox[0pt]{\mbox{\normalfont\tiny $P_{\theta}$}}}{\longrightarrow}}  b_I^{(3)} (\theta) = {\mathbb E_{\theta}} \Psi_3(X,Y,Z), \, n \to \infty.
$$

In efficiency calculations we shall not go beyond the location alternative other than few remarks on common parametric alternatives.  Let $P_{\theta}$ denote the alternative df $F(x, \theta) = F(x -\theta)$ with some symmetric df $F.$ Under these notations we obtain
\begin{multline*}
b_I^{(3)} (\theta)= P_{\theta}\{ |\min(Y,Z)| < |X|\} - P_{\theta}\{ |\max(Y,Z)|< |X|\} \\ = \int_{0}^{\infty}\left((1-F(-x -\theta))^2 - (1- F(x -\theta))^2 \right) d(F (x-\theta) -  F(-x -\theta))  \\- \int_{0}^{\infty}\left(F^2 (x-\theta) -  F^2 (-x -\theta) \right) d(F (x-\theta) -  F(-x -\theta)) \\ = 2\int_{0}^{\infty}\left(F(x-\theta) -  F(-x -\theta) \right)\left( 1 - F (x-\theta) -  F(-x -\theta)\right) d(F (x-\theta) -  F(-x -\theta)).
\end{multline*}
Assuming that $F$ is differentiable with the density $f$, we have for any $x$ and $\theta \to 0$
$$
\begin{array}{ll}
&F(x-\theta) -  F(-x -\theta) = 2F(x) - 1  + O(\theta^2), \\
&1 - F (x-\theta) -  F(-x -\theta) = 1 - F (x) -  F(-x) + 2\theta f(x) + O(\theta^2) \\&= 2f(x)\theta +  O(\theta^2).
\end{array}
$$
Consequently, under weak regularity conditions imposed on $F,$ we have
$$
b_I^{(3)} (\theta) \sim 8 \int_{0}^{\infty} (2 F(x) -1) f^2 (x) dx \cdot \theta, \, \theta \to 0.
$$
It follows that the local exact Bahadur slope \cite[Sect.7]{bah}, \cite{nik} is equivalent as $\theta \to 0$ to
$$
c_I^{(3)} (\theta) \sim 320 \left(\int_{0}^{\infty} (2 F(x) -1) f^2 (x) dx \right)^2 \cdot \theta^2.
$$

This local exact slope is equivalent to that of Litvinova's test studied  in \cite{lit}, \cite{lit1}. Her test was based on the Baringhaus-Henze characterization, and the test statistic appeared as $U$-statistic with the centered kernel
$$
\Phi(x,y,z) = \frac12 - \frac13\left({\bf{1}}\{ |\max(x,y)|< |z|\} + {\bf{1}}\{ |\max(x,z)|< |y|\} + {\bf{1}}\{ |\max(y,z)|< |x|\}  \right).
$$
The calculations are similar. While the limiting distributions have distinct variances and hence large deviation asymptotics are also different, the local exact slope is the same. Hence both tests are statistically equivalent for large samples, at least from the point of view of local Bahadur efficiency
(and also limiting Pitman efficiency).

According to the inequality (\ref{KL}), in our case of location
parameter we have under mild regularity conditions, see \cite{bah}, \cite{nik}, \ \S 4.4,
\begin{equation}
\label{BahRag}
320 \left(\int_{0}^{\infty} (2 F(x) -1) f^2 (x) dx \right)^2   \leq  I(f),
\end{equation}
where
$$
I(f)=\int_{-\infty}^{\infty}{\frac{\left(f'(x)\right)^2}{f(x)}dx}
$$
is the Fisher information.The local Bahadur  efficiency is equal to the ratio of the left and right parts in (\ref{BahRag}). Litvinova found rather high values of this efficiency for some concrete distributions. For instance, she found an efficiency as 0.977 for the normal distribution and 0.938 for logistic distribution. At the same time this efficiency is only 0.488 for the Cauchy distribution.

Close local  efficiencies appear for the skew alternative with the density $2f(x)F(\theta x),$ see \cite{azz}.  Litvinova \cite{lit1} explored also the contamination alternative and the Lehmann alternative, and the local efficiency has been always sufficiently high.

It follows that our test is also quite efficient with respect to the named alternatives.
Which test is better is to ascertain and can be explored either by power simulation or by the calculation of variances for corresponding $P-$values in the spirit of paper \cite{lam}.

\section{Kolmogorov-type test of symmetry for $k=2$}

\noindent In this section we consider the test of supremum type based on the test statistic
\begin{equation}
\label{Kolm}
D_n^{(2)} = \sup_t | G_n(t) - H_n(t) |.
\end{equation}
As this statistic is also distribution-free under $H_0$, we may assume  that the rv's $X_i$ are again uniformly distributed on $[-1,1]$ and that the supremum can be taken over $[-1,1].$

Limiting distribution and critical values of this  statistic are unknown but can be obtained via simulation. Therefore we will focus on its large deviations. Statistic (\ref{Kolm}) is the supremum of a {\it family of $U$-statistics} with the kernel of degree 2 depending on $t$, namely
\begin{equation}
\label{kernel}
   \Xi_2 (X,Y,t) = {\bf{1}}\{ |\min(X,Y)| < t\} - {\bf{1}}\{ |\max(X, Y)| < t\}, \, 0 \leq t \leq 1.
\end{equation}

In the following we need the projection function of the family of kernels (\ref{kernel}), see \cite{nikolm},
$$
\xi_2 (z,t)  = {\mathbb E} [\Xi_2 (X,Y,t) |Y=z] = {\mathbb P}\{ |\min(z, X)| < t\} - {\mathbb P}\{ |\max(z, X)| < t\}.
$$

This function depends on the relationship between $z$ and $t,$ and after some calculations we get
$$
\xi_2 (z; t) = \left\{\begin{array}{cc} - t, & {-1 \le z <-t,} \\
0, & {-t \le z \le t,} \\ t, & {t < z \le 1.}
\end{array} \right.
$$

Therefore we can calculate the so-called {\it variance function} \cite{nikolm} of the family of kernels (\ref{kernel}). We get
$$
\xi_2(t): = {\mathbb E} \xi_2^2(Y; t) = t^2 (1-t), \ 0 \le t \le 1.
$$
The maximum of this function is attained for $t= \frac23$ and is equal to $\frac{4}{27}$. We note that the variance function is non-degenerate in the sense of \cite{nikolm}, and hence we get, due to \cite{nikolm},  the  large deviation asymptotics
$$
\lim_{n \to \infty} n^{-1}\ln {\mathbb P}( D_n^{(2)}>a)  = - h_2(a)\sim
             - \frac{27}{32} a^2, \, \mbox{as} \, \, a \ \to 0,
$$
where $h_2$ is some analytic function in the neighbourhood of zero.

Hence the exact slope of our statistic $D_n^{(2)}$ \ is \ $2h_2(b_{D}^{(2)}(\theta))$, where
$$b_{D}^{(2)}(\theta) = \lim_{n \to \infty} D_n^{(2)}$$
a.s. under the alternative. Under the location alternative we can use the calculations made above and we get under minimal regularity assumptions
\begin{multline*}
b_{D}^{(2)} (\theta) = \sup_t | P_{\theta} \{ |\min(X,Y)| < t\} - P_{\theta} \{ |\max(X,Y)| < t\} \\
=\sup_t | (1-F(-t-\theta))^2 -(1-F(t-\theta))^2 - F^2(t-\theta) +F^2(-t - \theta)|  \\
=2\sup_t |(F(t-\theta) - F(-t - \theta))(1 - F(t-\theta) -F(-t - \theta))|  \\
\sim 4\sup_t |(2F(t) -1)|f(t)\cdot \theta, \, \theta \to 0.
\end{multline*}

Thus the local exact slope of the sequence $D_n^{(2)}$ satisfies the relation
$$
c_{D}^{(2)} (\theta) \sim 27 \sup_t (2F(t) -1)^2 f^2 (t) \cdot \theta^2, \, \theta \to 0.
$$
This local exact slope coincides with that of Kolmogorov-type  test from  \cite{bih} \ as evaluated in \cite{nikMMS}. The latter test is formally different being based on the difference
of $U$-empirical df's $F_n$ and $H_n$ but turns out to be asymptotically equivalent to our statistic $D_n^{(2)}$.

In any case, in \cite{nikMMS} the local Bahadur efficiency of both tests is calculated for location alternatives. It is 0.764 for the normal law,
0.750 for the logistic case, and 0.376 for the Cauchy distribution. For the Kolmogorov-type tests it is an adoptable result as such tests usually are less efficient than integral ones \cite{nik}.

\section{Integral tests in the general case}

\noindent We see that the used characterization of symmetry for $k=2$ leads to the tests which are asymptotically equivalent and equally efficient to known ones. Let us consider the general case when the tests are built on the characterization by the property
\begin{equation}
\label{Ahsan}
|\min(X_1,...,X_{k})| \mathrel{\stackrel{\makebox[0pt]{\mbox{\normalfont\tiny d}}}{=}} |\max(X_1,...,X_{k})|, \, k\ge 3.
\end{equation}

In the sequel the index $k$ or $k+1$ corresponds again to the degree of the kernel of $U-$statistic.
As in previous sections we associate with the condition (\ref{Ahsan}) of order $k$ the $U-$statistic of degree $k+1$
$$
I_n^{(k+1)} = { n\choose k+1 }^{-1}\sum_{1\le i_1< ... < i_{k+1} \le n} \Psi_{k+1} (X_{i_1}, ... , X_{i_{k+1}}),
$$
where the kernel $\Psi_{k+1}$ of degree $k+1$ is given after symmetrization  by
\begin{multline*}
(k+1)\Psi_{k+1}(X_1,...,X_{k+1})  \\= {\bf{1}}\{ |\min(X_1,...,X_k)| < |X_{k+1}|\} +...+ {\bf{1}}\{ |\min(X_2,...,X_{k+1})| < |X_{1}|\} \\ - {\bf{1}}\{ |\max(X_1,...,X_k)| < |X_{k+1}|\} -...- {\bf{1}}\{ |\max(X_2,...,X_{k+1})| < |X_{1}|\}.
\end{multline*}
In section 2 we studied the special case of this kernel for $k=2.$

When calculating the projection $\psi_{k+1}$ of this kernel, we are first interested in
$$
{\mathbb P}(|\min(X_1,...,X_{k-1},s)| < t) - {\mathbb P}(|\max(X_1,...,X_{k-1},s)| < t).
$$
Reasoning as above, we have for $s>0$ and $t \in [0,1]$
 $$
\medskip
\qquad \qquad {\mathbb P}(|\min(X_1,...,X_{k-1},s)| < t)= \left\{ \begin{array} {lll} (1+t)^{k-1} /2^{k-1}, &\mbox {\rm if} \, s \leq t,\\
 (1+t)^{k-1} /2^{k-1} - (1-t)^{k-1} /2^{k-1}, \, &\mbox {\rm if} \, \, t< s \leq 1.
\end{array}\right.
$$
Therefore, integrating, we get for $s> 0$
$$
{\mathbb P}(|\min(X_1,...,X_{k-1},s)| < |Z|) = \left( 2^{k} -2 + (1-s)^{k}\right) / k 2^{k-1}.
$$
In the same manner for $s\leq 0$ we obtain
$$
{\mathbb P}(|\min(X_1,...,X_{k-1},s)| < |Z|) =  \left(2^{k} -(1-s)^{k}\right)/k 2^{k-1}.
$$

Quite analogously we find the probabilities related to the maximum, namely
$$
\medskip
\qquad \qquad {\mathbb P}(|\max(X_1,...,X_{k-1},s)| < |Z|)= \left\{ \begin{array} {lll}
\left(2^{k} - (1 +s)^{k}\right) / k 2^{k-1}, \, s>0,\\
 \left(2^{k} -2 + (1+s)^{k}\right) / k2^{k-1},  \, s \leq 0.
\end{array}\right.
$$
Taking together our calculations, we obtain the projection of our kernel as
$$
 \psi_{k+1} (s)= \left\{ \begin{array} {lll} & \frac{(1+s)^{k} +(1-s)^{k} - 2}{(k +1)2^{k-1}}, \, \mbox {\rm if} & s> 0
;\\ &\frac{2 - (1+s)^{k} -(1-s)^{k}}{(k +1)2^{k-1}},  \, \mbox {\rm if} &s \leq 0.
\end{array}\right.
$$

Now we can calculate the variance $\sigma_{k+1}^2  = {\mathbb E} \psi_{k+1}^2 (X_1).$ It is given for any $k \ge 2$ by
$$
\sigma_{k +1}^2  = \frac{1}{2^{2k-2} (k +1)^2}\int_0^1 \left( (1+s)^{k } + (1-s)^{k } -2 \right)^2 \ ds > 0.
$$
In Table 1 we give some values of this variance which apparently has no nice explicit form.

\begin{table}[!hhh]\centering
\caption{Some exact values of the variance $\sigma^2_{k+1}.$ }

\medskip

\begin{tabular}{|c|c|}
\hline
 $k$ & Variance $\sigma_{k+1}^2$ \\
\hline
$k=2$ & 1/45\\
$k=3$ & 9/320 \\
$k=4$ & 2843/126000\\
$k=5$ & 2335/145152 \\
$k=6$ & 421691/37669632 \\
\hline
\end{tabular}
\end{table}

       By Hoeffding's theorem, see (\ref{normal}), the limiting distribution of $\sqrt{n}I_n^{(k+1)}$ is $N(0, (k+1)^2 \sigma^2_{k+1})$.  The large deviation asymptotics under $H_0$, see section 2, is given by
$$
\lim_{n\to \infty} n^{-1} \ln P ( I_n^{(k+1)} >a) = - f_{k+1} (a),
$$
where the function $f_{k+1}, k \geq 2,$ is analytic for sufficiently small $a>0,$ and such that
$$
f_{k+1} (a) \sim \frac{a^2}{2(k+1)^2 \sigma_{k+1}^2}, \quad \mbox{as} \,
\, a \to 0.
$$

Thus the local exact slope of the sequence of statistics $I_n^{(k+1)} , k\geq 2,$ is equivalent to
$$
c_{I}^{(k+1)} (\theta) \sim (b_{I}^{(k+1)} (\theta) )^2 / (k+1)^2 \sigma_{k+1}^2, \, \mbox {as} \,\, \theta \to 0.
$$
We  see that
\begin{multline*}
b_I^{(k+1)} (\theta)= P_{\theta}\{ |\min(X_1,...,X_{k}| < |Z|\} - P_{\theta}\{ |\max(X_1,...,X_{k}|< |Z|\}  \\ = \int_{0}^{\infty}\left((1-F(-x -\theta))^{k} - (1- F(x -\theta))^{ k } \right) d(F (x-\theta) -  F(-x -\theta))  \\- \int_{0}^{\infty}\left(F^{k} (x-\theta) -  F^{k} (-x -\theta) \right) d(F (x-\theta) -  F(-x -\theta))  \\ \sim  4k\int_{0}^{\infty}\left(F^{k-1} (x) -  F^{k-1}(-x) \right)f^2(x) dx \cdot \theta.
\end{multline*}

Hence, the local exact slope of the statistic of order $k$ is equal to
\begin{equation}
\label{locslo}
 c_{I}^{(k+1)} (\theta) \sim    \frac{16k^2}{(k +1)^2 \sigma_{k +1}^2}\left(\int_{0}^{\infty}\left(F^{k-1} (x) -  F^{k-1}(-x) \right)f^2(x) dx\right)^2 \cdot \theta^2.
\end{equation}

It is somewhat surprising to see that for $k=3$ we get from (\ref{locslo}), as $\theta \to 0:$
$$
c_{I}^{(4)} (\theta)\sim  c_{I}^{(3)} (\theta) \sim 320 \left(\int_{0}^{\infty} (2 F(x) -1) f^2 (x) dx \right)^2 \cdot \theta^2.
$$

But this equivalence is not long. Already for $k=4$ we get the variance $\sigma_5^2  = 2843/126000, $ hence the local exact slope is equivalent as $\theta \to 0$ to the expression
$$
c_{I}^{(5)} (\theta) \sim \frac{1290240}{2843} \left(\int_{0}^{\infty}\left(F^3 (x) -  F^3(-x)) \right)f^2(x) dx\right)^2\theta^2,
$$
which is different from the case $k=2$ and $k=3$.

For instance, in case of logistic distribution we have as $\theta \to 0$
$$
c_{I}^{(5)} (\theta) \sim \frac{1290240}{2843} \left(\int_{0}^{\infty}\frac{(e^{3x}-1)e^{2x}}{(e^{x} +1)^7} \ dx \right)^2 \theta^2 = \frac{1290240}{2843}\cdot\left(\frac{5}{192}\right)^2\theta^2  \approx 0.308 \cdot \theta^2.
$$
As the Fisher information  in this case  is $\frac13,$ the efficiency of our test is 0.925. This is high value comparable with the value 0.938 in case of lower dimensions $k=2$ and $k=3.$

In the case of normal law we get
$$
c_{I}^{(5)}(\theta) \sim \frac{1290240}{2843\cdot 4\pi^2} \left(\int_{0}^{\infty} \left(\Phi^3(x) - \Phi^3(-x)\right) \exp(-x^2) \ dx \right)^2 \theta^2 \approx 0.975 \cdot \theta^2.
$$
Note that 0.975 is just the value of local efficiency as the Fisher information is equal to 1. This is also high value. On the contrary, in the Cauchy case we get again much lower value of local efficiency 0.332.

It is interesting to compare the calculations of efficiencies for other common symmetric distributions and for other alternatives.

\section{ Local efficiency of \ Kolmogorov-type test in the \ general case}

\noindent Using the condition (\ref{Ahsan}) for any $k>2$, we can construct the Kolmogorov-type statistic $D_n^{(k)}$
according to (\ref{Kolm}). We concentrate here on large deviations and local efficiencies of such statistics for location alternatives.
It is necessary to consider the family of kernels, depending on $t \in [0,1]$ in a following way:
$$
\Psi_k (X_1,...,X_{k} ,t) =
{\mathbb P}(|\min(X_1,..., X_{k})| < t) - {\mathbb P}( |\max(X_1,..., X_{k})| < t).
$$

Let us calculate the projection of this family. We have
\begin{multline*}
\xi_k (z,t):= {\mathbb E}\left( \Psi(X_1,...,X_{k},t) | X_{ k}=z\right) \\={\mathbb P}(|\min(X_1,...,X_{k-1}, z)| < t) - {\mathbb P}( |\max(X_1,...,X_{k-1},z)| < t).
\end{multline*}

Using the calculations performed above,  we obtain
$$ \xi_k (z, t) = \left\{\begin{array}{cc}(1-t)^{k-1} /2^{k-1} - (1+t)^{k-1} /2^{k-1}, & {-1 \le z <-t,} \\
0, & {-t \le z \le t,} \\ (1+t)^{k-1} /2^{k-1} - (1-t)^{k-1} /2^{k-1}, & {t < z \le 1.}
\end{array} \right.$$

 Consequently the variance function is equal to
 \begin{multline*}
 \xi_k (t)= {\mathbb E} (\xi_k (Z,t))^2 =  \frac12\int_{-1}^{-t} ((1-t)^{k-1} /2^{k-1} - (1+t)^{k-1} /2^{k-1})^2 dx  \\ + \frac12\int_t^1 ( (1-t)^{k-1} /2^{k-1} - (1+t)^{k-1} /2^{k-1})^2 dx  \\= (1-t)((1+t)^{k-1} /2^{k-1} - (1-t)^{k-1} /2^{k-1})^2.
 \end{multline*}

 For $k=3$ we have again, as in the case $k=2$, the variance function
 $$\xi_3(t) = (1-t)t^2, -1 \le t \le 1 ,
 $$
 with the same maximum $\frac{4}{27},$ so that the large deviation asymptotics, see \cite{nikolm}, is given by the formula
 $$
 \lim_{n \to \infty} n^{-1}\ln {\mathbb P}( D_n^{(3)}> a)  =
            -h_3(a)= - \frac{3}{8} a^2 (1+o(1)),\ \mbox{as} \, \ a \to 0,
 $$
where $h_3$ is some analytic function in the vicinity of zero.

It is easy to see that the a.s. limit under the alternative of statistics $D_n^{(k)}$ admits the representation
$$
b_{D}^{(k)} (\theta) \sim 2k \sup_x f(x) [F^{k-1} (x) - F^{k-1} (-x) ]\cdot \theta, \, \theta \to 0.
$$
It follows that for $k=3$ the local exact slope has the form
$$
c_D^{(3)}(\theta) \sim 27\left(\sup_x [f(x)(2F(x) -1)]\right)^2 \cdot \theta^2, \, \ \theta \to 0,
$$
and the test is again equivalent to that of the case $k=2$ as in the instance of integral tests.

But in the case $k=4$ the situation changes as the variance function is
$$
\xi_4(t) = \frac{1}{16}(1-t)(3t + t^3)^2, \, \, 0 \le t \le 1.
$$
We find numerically that the maximum of the variance function is equal to 0.1123... \. Hence the large deviation result is different and reads
$$
\lim_{n \to \infty} n^{-1}\ln  {\mathbb P}( D_n^{(4)} > a)  =  - h_4(a) = - 0.2783... \cdot a^2 (1+o(1)), \, \mbox{as} \, a \to 0.
$$
Therefore the  exact slope admits the representation
$$
c_{D}^{(4)}(\theta) \sim 35.622... \sup_x \left[ f (x) \left(F^3 (x) - F^3 (-x)\right) \right]^2 \cdot \theta^2, \, \theta \to 0.
$$

In case of logistic distribution and $k=4$ we have in the right-hand side
$$
35.622... \sup_x \left(\frac{e^{x}(e^{3x}-1)}{(1+e^{x})^5}\right)^2
 \approx 0.232,
$$
that gives for local efficiency lower result 0.696 than in previous cases.

For the normal law we find that
$$
\frac{1}{2\pi}\sup_x \ e^{-x^2} \left(\Phi^3(x) -\Phi^3(-x)\right)^2 \approx 0.0206.
$$
Consequently, the efficiency is approximately $0.733.$ Similar calculations show that for the Cauchy law the local efficiency equals 0.313. All these efficiencies are reasonable but moderate.

\section{Discussion}

\noindent We can resume the calculations of efficiencies in table 2. One sees  that for logistic and normal distributions the values of efficiencies of integral tests for location alternative are rather high in comparison with other nonparametric tests of symmetry, see \cite[Ch.4]{nik}.

\begin{table}[h]
\caption{Local Bahadur efficiencies in location case.}\label{lae}
\bigskip
 \hbox to
\textwidth{\hfill\begin{tabular}{|l|c|c|c|} \hline & & &  \\ Statistic/Density
 & Logistic  & Normal   & Cauchy  \\

  \hline & & & \\ $I_n^{(3)},I_n^{(4)} $ & 0,938 & 0,977 & 0,488  \\ \hline \
&\ &\ &\\ $I_n^{(5)}$ & 0,925 & 0,975 & 0,332 \\ \hline \ &\ &\
&\\ $D_n^{(2)},D_n^{(3)}$ & 0,750 & 0,764 & 0,376  \\ \hline \ &\ &\ & \\
$D_n^{(4)}$ & 0,696 & 0,733 & 0,313  \\ \hline

\end{tabular}\hfill}
\end{table}

At the same time the results for the Cauchy law are mediocre. It would be of interest to study other alternatives and to compare the efficiency values with the power simulations for moderate sample size.

The efficiencies of Kolmogorov-type tests are lower but have tolerable values. One should keep in mind that these
tests are always consistent while the integral tests of structure (\ref{integral}) have mostly one-sided character, and their consistency depends on the alternative.

We can also presume the deterioration of efficiency properties for our tests with the growth of their order and degree of complexity, at least for location alternative. Hence the simplest test statistics $I_n^{(3)}$ and $D_n^{(2)}$ and their equivalents described above seem to be most suitable for practical use.


\begin{thebibliography}{99}

\bibitem{ahsan} Ahsanullah, M.: On some characteristic property of symmetric distributions.
Pakist. J. Statist. {\bf 8}, 19 -- 22(1992).

\bibitem{azz} Azzalini, A. with the collaboration of Capitanio, A.: The skew-normal and related fa\-milies. Cambridge University Press, New York (2014).

\bibitem{Bah67} Bahadur,~R.~R.: Rates of convergence
of estimates and test statistics. Ann.\ Math.\ Statist., {\bf 38}, 303 -- 324(1967).

\bibitem{bah} Bahadur, R.R.: Some limit theorems in statistics. SIAM, Philadelphia(1971).

\bibitem{bih}  Baringhaus, L., Henze, N.: A characterization of
and new consistent tests for symmetry. Commun. Statist.- Theory
Meth., \ {\bf 21}, 1555 -- 1566(1992).

\bibitem{halm} Halmos, P. R.: The theory of unbiased estimation. Ann. Math.
Statist., {\bf 17}, 34 -- 43(1946).

\bibitem{helm} Helmers, R., Janssen, P., Serfling, R.: Glivenko-Cantelli properties of some generalized
empirical df's and strong convergence of generalized $L$-statistics. Probab. Theor. Rel.
Fields, {\bf 79}, 75 -– 93(1988).

\bibitem{henMei}
 Henze, N., Meintanis, S.: Goodness-of-fit tests based on a
 new characterization of the exponential distribution. Commun. Statist. Theor. Meth.,
 {\bf  31}, 1479 -- 1497(2002).

\bibitem{hoe}
 Hoeffding, W.: A class of statistics with asymptotically normal
distribution.  Ann.\ Math.\ Statist., {\bf 19}, 293 -- 325(1948).

\bibitem{kor}  Korolyuk, V.S., Borovskikh, Yu.V.: Theory of  $U$-statistics.
 Kluwer, Dordrecht(1994).

 \bibitem{lam} Lambert, D., Hall, W.J.: Asymptotic Lognormality of $P-$Values. Ann. Stat.
{\bf 10}, 44 -- 64(1982).

\bibitem{Lee} Lee~A.~J.: $U-$statistics: Theory and Practice. \ Dekker, New York(1990).


\bibitem{Lin} Linnik, Yu. V.: Linear forms and statistical criteria. \ I,\ II.
Ukrainian Math. J.  {\bf 5}, 207 -- 243; {\bf 5}, 247 -- 290(1953).

\bibitem{lit}Litvinova, V. V.:  New nonparametric test for symmetry and its asymptotic efficiency. Vestnik of St.Petersburg University. Mathematics. {\bf 34},  12 -- 14(2001).

\bibitem{lit1} Litvinova, V.V.: Asymptotic properties of goodness-of-fit and symmetry tests based on characterizations. Ph.D. thesis. Saint-Petersburg University (2004).

\bibitem{morSz}  Morris, K., Szynal, D.: Goodness-of-fit tests using
characterizations of continuous distributions.  Appl. Math. (Warsaw). {\bf 28}, 151 -- 168(2001).

\bibitem{mulnik} Muliere, P., Nikitin, Ya. Yu.: Scale-invariant test of normality
based on Polya's characterization. Metron. {\bf 60},  21 -- 33(2002).

\bibitem{nik} Nikitin, Y.:  Asymptotic efficiency of nonparametric tests.
Cambridge University Press, New York(1995.)

\bibitem{nikMMS} Nikitin, Ya. Yu.: On Baringhaus-Henze test for symmetry: Bahadur efficiency and local optimality for shift alternatives. Math. Methods of Statist. {\bf 5}, 214 -- 226(1996).

\bibitem{nik96} Nikitin, Ya. Yu.:  Bahadur efficiency of a test of exponentiality based on a loss-of-memory
type functional equation. J.\ Nonparam.\ Statist. {\bf 6}, 13 -- 26(1996).

\bibitem{nikolm} Nikitin, Ya. Yu.: Large deviations of $U$-empirical
Kolmogorov-Smirnov tests, and their efficiency.  J.\ Nonparam. \ Statist., {\bf 22
}(2010), 649--668(2010).

\bibitem{nikiponi} Nikitin, Ya. Yu., Ponikarov, E.V.: Rough large deviation asymptotics of Chernoff type for von Mises functionals and $U$-statistics.   Proc. of St.Petersburg Math. Soc.  {\bf 7}, 124 -- 167(1999). Engl. transl. in AMS Transl., ser.2., {\bf 203}, 107 -- 146(2001).

\bibitem{nikVol} Nikitin, Ya. Yu., Volkova, K.Yu.:  Asymptotic efficiency of
exponentiality tests based on order statistics characterization.  Georgian Math. J., {\bf 17}, 749 -- 763(2010).




\end{thebibliography}
\end{document}